\documentclass[11pt]{article}
\usepackage{amsmath,amsfonts,amsthm,amssymb}
\usepackage{graphicx,url,subfig}
\usepackage[colorlinks=true,citecolor=blue]{hyperref}
\usepackage{enumerate}
\graphicspath{ {./Images/} }
\usepackage{color,tikz,mathtools}
\usepackage{float,wrapfig,fullpage}
\usepackage{enumerate, enumitem}

\usepackage{csquotes}

\usepackage{thm-restate}

\usepackage[margin = 2.50cm]{geometry}

\newcommand{\R}{\mathbb{R}}
\newcommand{\N}{\mathbb{N}}
\newcommand{\F}{\mathcal{F}}

\newcommand{\B}{\mathcal{B}}

\newtheorem{theorem}{Theorem}[section]
\newtheorem{corollary}[theorem]{Corollary}

\newtheorem{lemma}[theorem]{Lemma}

\newtheoremstyle{case}{}{}{}{}{}{:}{ }{}
\theoremstyle{case}
\newtheorem{case}{Case}

\newcommand{\NN}{\mathbb{N}}

\newcommand{\QQ}{\mathbb{Q}}
\newcommand{\RR}{\mathbb{R}}

\newcommand{\cB}{\mathcal{B}}

\newcommand{\cF}{\mathcal{F}}

\newcommand{\cH}{\mathcal{H}}

\newcommand{\cK}{\mathcal{K}}

\newcommand{\cS}{\mathcal{S}}
\newcommand{\cT}{\mathcal{T}}

\newcommand{\remove}[1]{}

\title{
    Stabbing boxes with finitely
    many axis-parallel lines and flats
}

\author{
Sutanoya Chakraborty\footnote{Indian Statistical Institute, India}
\and
Arijit Ghosh\footnotemark[1]
\and
Soumi Nandi\footnotemark[1]
}

\date{}

\begin{document}

\maketitle

\begin{abstract}
    In this short note, we provide the necessary and sufficient condition for an infinite collection of axis-parallel boxes in $\R^{d}$ to be pierceable by finitely many axis-parallel $k$-flats, where $0 \leq k < d$. We also consider {\em colorful} generalizations of the above result and establish their feasibility. The problem considered in this paper is an infinite variant of the Hadwiger-Debrunner $(p, q)$-problem. 
\end{abstract}

\section{Introduction}
\label{sec:introduction}

Let $\cF$ be a collection of subsets of $\R^{d}$. We say a subset $T$ of $\R^{d}$ {\em pierces} (or \emph{hits} or {\em stabs}) $\cF$ if for all $S \in \cF$ we have $S\cap T \neq \emptyset$. A family $\cT$ of subsets of $\R^{d}$ is said to be a \emph{transversal} for a family $\cF$ of sets in $\RR^d$ if for each $F\in\cF$ there exists a $\tau\in\cT$ such that $F\cap\tau\neq\emptyset$. Additionally, we say $\mathcal{T}$ is a {\em finite-size transversal} of $\cF$ if $\mathcal{T}$ is a transversal of $\cF$ and the size of $\mathcal{T}$ is finite.
$\cT$ can be a family of points or lines or hyperplanes or $k$-flats\footnote{By $k$-flats we mean $k$-dimensional affine subspaces of $\R^{d}$.} etc. If $\mathcal{T}$ is a collection of $k$-flats (or axis-parallel $k$-flats) in $\RR^{d}$ then $\mathcal{T}$ is called a {\em $k$-transversal} (or axis-parallel $k$-transversal) of $\cF$, and additionally, if $k$ is $0$ or $d-1$, then $\mathcal{T}$ is called a point transversal or a hyperplane transversal of $\cF$ respectively.

Helly's Theorem~\cite{Helly23}, a fundamental result in discrete and convex geometry, states that if $\cF$ is a family of compact convex sets in $\RR^d$ such that every $d+1$ sets from $\cF$ can be pierced by a single point, then the whole family $\cF$ is pierceable by a single point. The compactness requirement can be dropped if $\cF$ is finite. Lov\'asz proved the following {\em colorful} variant of the Helly Theorem, though the proof first appeared in a paper by B\'ar\'any~\cite{Barany82}. 
\begin{theorem}[B\'ar\'any~\cite{Barany82}:
Colorful Helly Theorem]
Suppose $\cF_1,\dots,\cF_{d+1}$ are non-empty finite families of convex sets in $\RR^d$ such that for all $\left( C_{1}, \dots, C_{d+1}\right) \in \cF_{1} \times \dots \times \cF_{d+1}$ we have $\bigcap_{i \in [d+1]} C_{i} \neq \emptyset$. Then there exists an $i\in [d+1]$ and a point in $\R^{d}$ that pierces all the sets in $\cF_i$.   
\end{theorem}
\noindent
Over the last hundred years, several generalizations and applications of Helly's Theorem have been extensively studied. For an overview of the same, one can see Amenta et al.~\cite{AmentaDeLS17}, De Loera et al.~\cite{de2019discrete}, and B{\'a}r{\'a}ny and Kalai~\cite{BaranyK2022helly}.

Generalizing Helly's Theorem to general $k$-transversals, for $1\leq k \leq d-1$, has been an active area of research. See the survey by Holmsen and Wenger~\cite{HolmsenW17} for more details. 
Santal\'{o}~\cite{Santal2009} showed that there can be no {\em Helly-type} theorem for arbitrary families of compact convex sets in $\R^3$ and line transversals. 
Because of Santal\'{o}'s impossibility result, the follow-up works on 
$k$-transversals are for families of convex sets that satisfy additional properties~\cite{Hadwiger56,Hadwiger57,DGK,GoodmanPJAMS88,ABJP,AGP,HolmsenKL03,CheongGHP08}.

For $p \geq q$, we say a family $\cF$ of convex sets in $\R^{d}$ satisfies the $(p,q)$-property for points if every $p$ sets in $\cF$ contain $q$ sets that can be pierced by a point. Hadwiger and Debrunner~\cite{HadwigerD57} proved the following fundamental generalization of Helly's Theorem.
\begin{theorem}[Hadwiger and Debrunner~\cite{HadwigerD57}] Let $p \geq q \geq d + 1$ with $q > \frac{(d-1)p + d}{d}$, and $\cF$ be a finite family of compact convex sets in $\R^{d}$ satisfying the $(p, q)$-property for points.
Then $\cF$ can be pierced by $p - q + 1$ points.
\label{theorem-Hadwiger-Debrunner-problem}
\end{theorem}
\noindent
Hadwiger and Debrunner conjectured that for all $p \geq q \geq d + 1$, if a family $\cF$ of compact convex sets in $\R^{d}$ satisfies the $(p, q)$-property, then $\cF$ has a point transversal of constant size. 
This conjecture was resolved by Alon and Kleitman~\cite{AlonK92} $35$ years later. 

\begin{theorem}[Alon and Kleitman~\cite{AlonK92}: $(p,q)$-Theorem]
    For any three natural numbers $p\geq q\geq d+1$, $\exists c=c(p,q,d)$ such that if $\cF$ is a collection of compact convex sets in $\RR^d$ satisfying the $(p,q)$-property for points, then there exists a point transversal for $\cF$ with size at most $c$.
\end{theorem}
\noindent
We can similarly define $(p,q)$-property for $k$-flats. For $p \geq q$, a family $\cF$ of convex sets in $\R^{d}$ satisfies the $(p,q)$-property for $k$-flats if every $p$ sets in $\cF$ contain $q$ sets that can pierced by a single $k$-flat.
Alon and Kalai~\cite{AlonK95} proved the $(p,q)$-theorem for hyperplane transversals.
\begin{theorem}[Alon and Kalai~\cite{AlonK95}: $(p,q)$-Theorem for hyperplane transversals]
    \label{thm-Alon-Kalai-pq-hyperplane-transversal}
    For any three natural numbers $p\geq q\geq d+1$, $\exists c'=c'(p,q,d)$ such that if $\cF$ is a collection of compact convex sets in $\RR^d$ satisfying the $(p,q)$-property for hyperplanes then there exists a hyperplane transversal for $\cF$ with size at most $c'$.
\end{theorem}
\noindent
Alon, Kalai, Matou\v{s}ek, and Meshulam~\cite{AKMM} proved the impossibility of getting a $(p,q)$-theorem for $k$-transversal when $1\leq k < d-1$. Since its introduction by Hadwiger and Debrunner~\cite{HadwigerD57}, this relaxed version of Helly's Theorem has been an active area of research in discrete and convex geometry, and is now popularly known as {\em Hadwiger-Debrunner $(p,q)$-problem} or just {\em $(p,q)$-problem}.

A family $\cF$ is said to satisfy the $\left(\aleph_{0}, q\right)$-property, where $q: = q(k,d) \geq k+2$, for $k$-transversal if for any infinite sequence 
$\left\{ C_{i}\right\}_{i \in \mathbb{N}}$ of sets from $\cF$ there exist $q$ sets from the sequence that can be pierced by a single $k$-flat. 
Keller and Perles~\cite{KellerP22} were the first to introduce this new infinite variant of the $(p,q)$-property and proved the following results:

\begin{theorem}[Keller and Perles~\cite{KellerP22}: $(\aleph_{0},k+2)$-Theorem for $k$-transversals]\label{Chaya}
   Let $\cF$ be a family of compact convex sets in $\RR^d$ such that each $C\in\cF$ contains a ball of radius $r>0$ and is contained in a ball of radius $R$, and $0\leq k<d$. If $\cF$ satisfies the $(\aleph_{0},k+2)$-property for $k$-flats, then $\cF$ has a finite-size $k$-transversal.
\end{theorem}

\begin{theorem}[Keller and Perles~\cite{KellerP22}: $(\aleph_{0},2)$-Theorem for point transversals]\label{th:H_balls_pt}
    Suppose $\cF$ is a family of closed balls in $\R^d$ (with no restriction on the radii) satisfying the $\left( \aleph_{0}, 2\right)$-property for points. 
    Then $\cF$ has a finite-size point transversal.
\end{theorem}
\noindent
In a later work, Keller and Perles~\cite{KellerP23} introduced the notion of {\em near-balls}. A collection $\cF$ of subsets of $\R^{d}$ is called a {\em near-ball collection} if there exists a constant $\alpha > 0$ such that for all $C \in \cF$, there exist closed balls $B(p_{C}, r_{C})$ and $B(p_{C}, R_{C})$ satisfying the following two conditions:
\begin{eqnarray}
    &B(p_{C}, r_{C}) \subseteq C \subseteq B(p_{C}, R_{C}),~\mbox{and}&  \\
    &R_{C} \leq \min \left\{ r_{C} + \alpha, \alpha r_{C}\right\}.&
\end{eqnarray}
where $B(p, r)$ denotes a closed ball in $\R^d$ with center at $p$ and radius $r$.
Keller and Perles~\cite{KellerP23} proved the following result:
\begin{theorem}[Keller and Perles~\cite{KellerP23}: $(\aleph_{0},k+2)$-Theorem for near balls]\label{Keller-Perles-nearballs-2023}
    Let $\cF$ be a near-ball collection of compact sets from $\cF$ satisfying the $(\aleph_{0}, k+2)$-property for $k$-flats. Then $\cF$ has a finite-size $k$-transversal. 
\end{theorem}
\noindent
Clearly, Theorem~\ref{Chaya} directly follows from Theorem~\ref{Keller-Perles-nearballs-2023}, and the following corollary of Theorem~\ref{Keller-Perles-nearballs-2023} is a generalization of Theorem~\ref{th:H_balls_pt} to general $k$-flats.  

\begin{corollary}[Keller and Perles~\cite{KellerP23}: $(\aleph_{0},k+2)$-Theorem for closed balls]\label{cor:H_balls_pt_k_flats}
    Suppose $\cF$ is a family of closed balls in $\R^d$ (with no restriction on the radii) satisfying the $\left( \aleph_{0}, k+2\right)$-property for $k$-flats. 
    Then $\cF$ has a finite-size $k$-transversal.
\end{corollary}
\noindent
Jung and P\'{a}lv\"{o}lgyi~\cite{JungP2023} 
recently\footnote{Jung and P\'{a}lv\"{o}lgyi~\cite{JungP2023} uploaded their manuscript to arXiv a few months after the first draft of this paper~\cite{chakraborty2023stabbing} was uploaded to arXiv. \remove{Need to say in the footnote that this paper appeared after our result.}} developed a new framework using which one can lift $(p, q)$-type theorems for finite families to prove new $(\aleph_{0}, q)$-Theorems for infinite families, and also provided new proofs of already existing results. Chakraborty, Ghosh and Nandi~\cite{chakraborty2024heterochromatic,chakraborty2024countably} have also given colorful generalizations of the above theorems. 

\paragraph{Our results:}
Studying geometric, combinatorial, algorithmic, and even topological properties of a collection of axis-parallel boxes has been an active area of research in discrete and computational geometry~\cite{ChazelleL01,Toth08,ChenPST09,PachT10,RafalinST10,DumitrescuT15,ChalermsookW21,Mitchell21,ChanH21,GalvezKMMPW22}. 
{\bf For the remainder of the paper, unless otherwise stated explicitly, the term ``boxes'' is used to mean axis-parallel boxes.}
In the context of transversal theory, boxes have been one of the fundamental objects of study~\cite{DanzerG82,GyarfasL85,Fon-Der-FlaassK93,Alon98,Nielsen00,KatzNS03,PachT12,Alon12,PachT13,CorreaFPS15,KupavskiiMP16,MustafaR17,MustafaDG18,ChudnovskySZ18,KellerS20,ChenD20,tomon2022lower}. We build on these previous works by studying the $(\aleph_0,q)$-problem for boxes and axis-parallel flats.

An infinite collection $\mathcal{G}$ of subsets of $\R^{d}$ is said to satisfy the {\em $(\aleph_{0},2)$-property for axis-parallel $k$-flats} if for every infinite sequence $\left\{ S_{n}\right\}_{n \in \mathbb{N}}$, where $S_{n} \in \mathcal{G}$ for all $n \in \mathbb{N}$, there exist $S_{i}$ and $S_{j}$ with $i\neq j$ from $\left\{ S_{n}\right\}_{n \in \mathbb{N}}$ that can be pierced by a single axis-parallel $k$-flat. 
In this short note, we prove the following $(\aleph_{0},2)$-Theorem for boxes.

\begin{restatable}[$(\aleph_{0},2)$-Theorem for boxes and 
axis-parallel $k$-flats]{theorem}{kflats}
Let $\cF$ be an infinite collection of boxes in $\R^{d}$ and $0 \leq k < d$.
If $\cF$ satisfies the $(\aleph_{0}, 2)$-property for axis-parallel $k$-flats, then 
$\cF$ has a finite-size axis-parallel $k$-transversal. 
\label{thm-main-theorem}\end{restatable}


%
\noindent
It can be easily seen that the axis-parallel condition in Theorem~\ref{thm-main-theorem} cannot be relaxed if we still want the conclusion to hold.
Jung and P\'{a}lv\"{o}lgyi~\cite{JungP2023} have shown the following result about convex sets and axis-parallel hyperplanes. 
\begin{theorem}[Jung and P\'{a}lv\"{o}lgyi~\cite{JungP2023}: Convex sets and hyperplane transversals]
    Let $\cF$ be a collection of compact convex sets in $\R^{d}$ that satisfies the $(\aleph_{0}, 2)$-property for axis-parallel hyperplanes. Then $\cF$ has a finite-size hyperplane transversal. 
    \label{thm:JungPconvexset-axis-parallel-hyperplane}
\end{theorem}

\noindent
From Theorem~\ref{thm-main-theorem} we directly get the following generalization of Theorem~\ref{thm:JungPconvexset-axis-parallel-hyperplane}.
\begin{corollary}\label{corollary-hyperplane-connected-sets}
    Let $\F$ be a collection of compact connected sets in $\R^d$ that satisfies the $(\aleph_0,2)$-property for axis-parallel hyperplanes. Then there exists a finite collection of axis-parallel hyperplanes that is a transversal for $\F$.
\end{corollary}
\noindent
For a proof of the above corollary, it suffices to observe that an axis-parallel hyperplane intersects a compact connected set $C$ in $\R^d$ if and only if it intersects the smallest axis-parallel box in $\R^d$ that contains $C$. 

\begin{figure}
        \centering
        \includegraphics[width=0.8\linewidth]{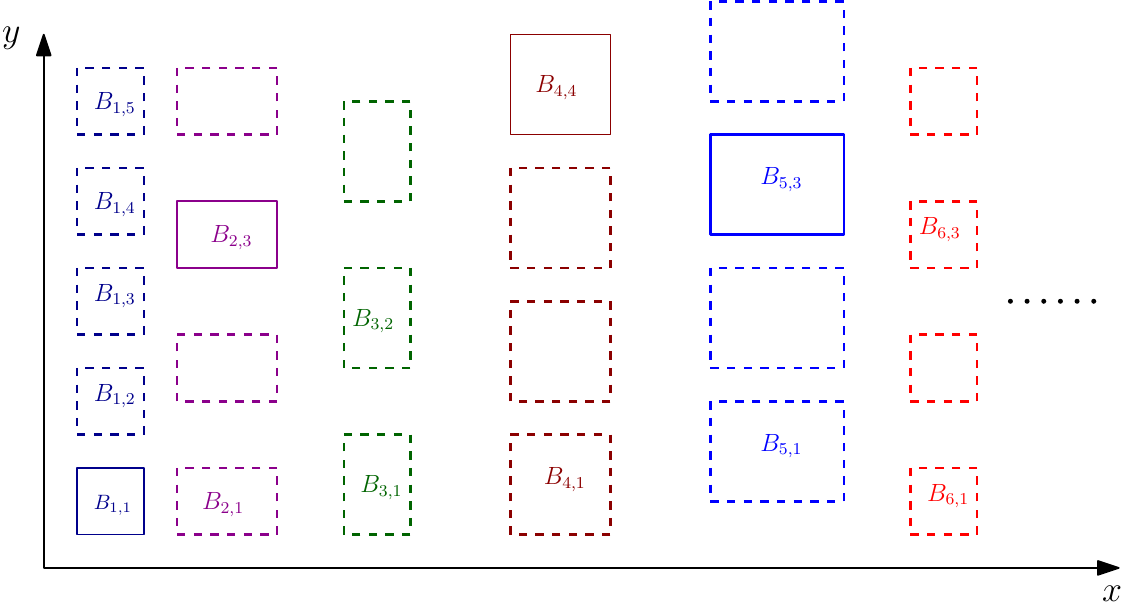}
        \caption{$\{B_{1,1}, B_{2,3},B_{4,4}, B_{5,3},\dots\}$ is a colorful sequence with respect to $\{\cF_n\}_{n\in\NN}$, where $\cF_n=\{B_{n,m}\;|\;m\in\NN\}$}
        \label{fig:enter-label3}
    \end{figure}

Given an infinite sequence $\left\{ \mathcal{G}_{n}\right\}_{n \in \mathbb{N}}$ of families of subsets of $\R^{d}$, we say that an infinite sequence $\left\{ S_{m}\right\}_{m \in \mathbb{N}}$ is {\em a colorful sequence} with respect to $\left\{ \mathcal{G}_{n}\right\}_{n \in \mathbb{N}}$ if the following two conditions hold:
    \begin{itemize}
        \item
            There is an infinite sequence $\{n_m\}_{m\in\N}$ in $\N$ with the property that $n_i>n_j$ if $i>j$, and
        \item 
            for all $m \in \N$, we have $S_{m} \in \mathcal{G}_{n_m}$.
%
%
    \end{itemize}
\noindent For example, in Figure~\ref{fig:enter-label3}, $\{B_{1,1}, B_{2,3},B_{4,4}, B_{5,3},\dots\}$ is a colorful sequence with respect to $\{\cF_n\}_{n\in\NN}$, where $\cF_n=\{B_{n,m}\;|\;m\in\NN\}$.
    
An infinite sequence $\left\{ \mathcal{G}_{n} \right\}_{n \in \mathbb{N}}$ is said to satisfy the {\em colorful $(\aleph_{0},2)$-property for axis-parallel $k$-flats} if for every colorful sequence $\left\{ S_{n}\right\}_{n \in \mathbb{N}}$ with respect to $\left\{ \mathcal{G}_{n} \right\}_{n \in \mathbb{N}}$, there exist two sets $S_{i}$ and $S_{j}$ with $i \neq j$ from $\left\{ S_{n}\right\}_{n \in \mathbb{N}}$ which can be pierced by a single axis-parallel $k$-flat.
We prove the following {\em colorful} generalization of Theorem~\ref{thm-main-theorem}. 
\begin{restatable}[Colorful generalization of Theorem~\ref{thm-main-theorem}]{theorem}{colorfulmaintheorem}
    Let $\left\{\cF_{n}\right\}_{n\in \NN}$ be an infinite sequence of families of axis-parallel boxes in $\RR^d$ and $0 \leq k < d$. If the infinite sequence $\left\{\cF_{n}\right\}_{n\in \NN}$ satisfies the colorful $(\aleph_{0},2)$-property for axis-parallel $k$-flats, then there exists an $\ell \in \mathbb{N}$ such that 
    $\cF_{\ell}$ has a finite-size axis-parallel $k$-transversal.
    \label{thm-colorful-main-theorem}
\end{restatable}
\noindent Using the observations of Chakraborty, Ghosh and Nandi~\cite[Section~3.1]{chakraborty2024heterochromatic}, we get the following stronger result. Note that completeness, we have reproduced the arguments of Chakraborty, Ghosh, and Nandi~\cite{chakraborty2024heterochromatic} in Appendix~\ref{sec:appendix1}.


\begin{corollary}
    Let $\{\F_n\}_{n\in\NN}$ be an infinite sequence of families of axis-parallel boxes in $\R^d$, and $0\leq k<d$.
    If $\{\F_n\}_{n\in\NN}$ satisfies the colorful $(\aleph_0,2)$-property for axis-parallel $k$-flats, then there exists an $N\in\NN$ and a finite collection $\cK$ of axis-parallel $k$-flats that is a transversal for all $\F_n$ with $n\geq N$.
\end{corollary}
\noindent
Putting $k=d-1$, the above corollary also provides a colorful generalization of Corollary~\ref{corollary-hyperplane-connected-sets}.
\begin{corollary}[Generalization of Corollary~\ref{corollary-hyperplane-connected-sets}]
   Let $\{\F_n\}_{n\in\NN}$ be an infinite sequence of families of compact connected sets in $\R^d$.
   If $\{\F_n\}_{n\in\NN}$ satisfies the colorful $(\aleph_0,2)$-property for axis-parallel hyperplanes, then there exists an $N\in\NN$ and a finite collection $\cK$ of axis-parallel hyperplanes that is a transversal for all $\F_n$ with $n\geq N$.
   \label{corollary-colorful-hyperplane-connected-sets}
\end{corollary}

\begin{figure}
    \centering
    \includegraphics[width=0.8\linewidth]{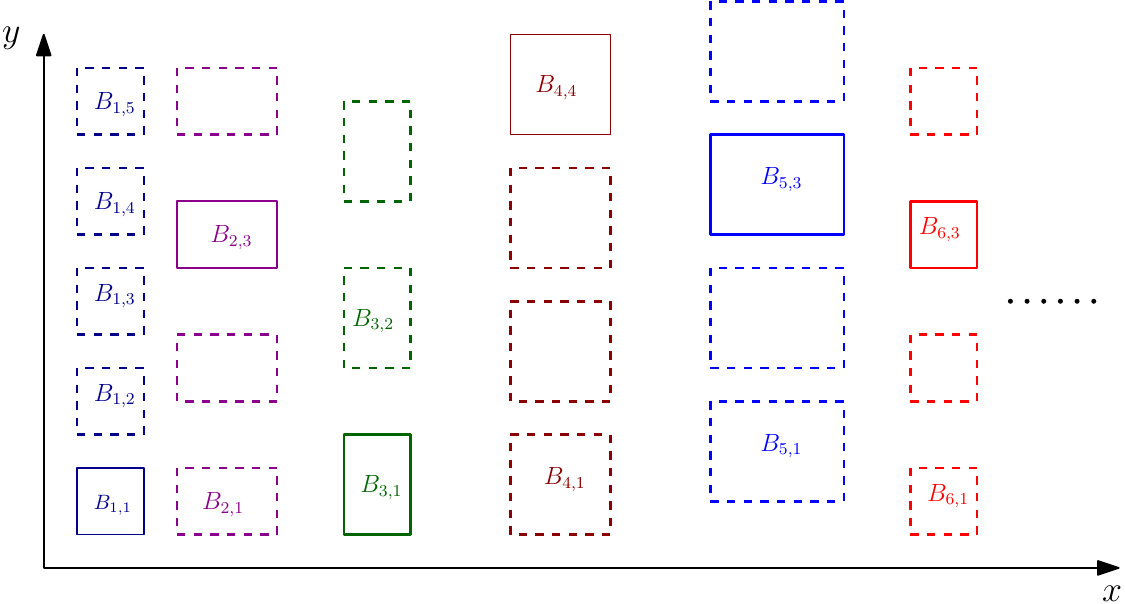}
    \caption{$\{B_{1,1}, B_{2,3}, B_{3,1}, B_{4,4}, B_{5,3}, B_{6,3}, \dots\}$ is a strongly colorful sequence with respect to $\{\cF_n\}_{n\in\NN}$, where $\cF_n=\{B_{n,m}\;|\;m\in\NN\}$}
    \label{fig:enter-label2}
\end{figure}

Given an infinite sequence $\left\{ \mathcal{G}_{n}\right\}_{n \in \mathbb{N}}$ of families of subsets of $\R^{d}$, we say that an infinite sequence of subsets $\left\{ S_{n}\right\}_{n \in \mathbb{N}}$ from $\R^{d}$ is a {\em strongly colorful sequence} of $\left\{ \mathcal{G}_{n}\right\}_{n \in \mathbb{N}}$ if for all $n \in \NN$, we have $S_{n} \in \mathcal{G}_{n}$. For example, in Figure~\ref{fig:enter-label2}, $\{B_{1,1}, B_{2,3}, B_{3,1}, B_{4,4}, B_{5,3}, B_{6,3}, \dots\}$ is a strongly colorful sequence with respect to $\{\cF_n\}_{n\in\NN}$, where $\cF_n=\{B_{n,m}\;|\;m\in\NN\}$

Given Theorem~\ref{thm-main-theorem} and Theorem~\ref{thm-colorful-main-theorem}, it is natural to ask if even a {\em stronger colorful} variant along the lines given below also holds:
\begin{center}
    \emph{
        Let $\left\{\cF_{n}\right\}_{n\in \NN}$ be an infinite sequence of families of axis-parallel boxes in $\RR^d$ and $0 \leq k < d$. If every strongly colorful sequence $\left\{ B_{n} \right\}_{n \in \NN}$ with respect to $\left\{ \cF_{n} \right\}_{n \in \NN}$ contain two boxes $B_{i}$ and $B_{i}$,
        with $i \neq j$, from $\left\{ B_{n} \right\}_{n \in \NN}$ and both $B_{i}$ and $B_{j}$ can be pierced by a single axis-parallel
        $k$-flat then there exists a $\ell \in \NN$ such that $\cF_{\ell}$ has a constant size axis-parallel $k$-transveral.
    }
\end{center}
\noindent
In the next theorem, we construct a sequence of families of boxes in $\mathbb{R}^{d}$ showing the infeasibility of the above statement.
\begin{restatable}[Impossibility of strong colorful generalization of Theorem~\ref{thm-main-theorem}]{theorem}{impossibility}
    \label{impossibility_strong_colorful}
    For all $d\in\N$ and $0 \leq k< d$, there exists an infinite sequence $\left\{\cF_{n}\right\}_{n \in \mathbb{N}}$ of families of boxes in $\RR^{d}$ satisfying the following properties:
    \begin{itemize}
        \item   
            for all $n \in \mathbb{N}$, $\cF_{n}$ does not have a finite size axis-parallel $k$-transversal, and
            
        \item
            for any $t\in\NN$, every infinite sequence $\left\{ B_{n}\right\}_{n \in \mathbb{N}}$, where $B_{n} \in \cF_{n+t}$ for all $n \in \mathbb{N}$, contains two boxes $B_{i}$ and $B_{j}$ with $i \neq j$ and both $B_{i}$ and $B_{j}$ can be pierced by a single axis-parallel $k$-flat. 
    \end{itemize}
\label{thm: impossibility}
\end{restatable}
\noindent
Keller has also given an independent construction ensuring the same impossibility as the above theorem, see~\cite{JungP2023}.

%

\paragraph{Notations:}
Throughout the rest of the paper, we will be using the following notations:
\begin{itemize}
    
    \item 
        For all $n \in \N$, $[n]$ denotes the set $\left\{ 1, \dots, n\right\}$.

    \item       
        $\mathbb{Q}$ is the set of rational numbers, and 
        $$
            \mathbb{Q}^{d} : = \underbrace{\mathbb{Q} \times \dots \times \mathbb{Q}}_{d~times}.
        $$
    
    \item 
        For all $a, b \in \R$ with $a \leq b$, $[a,b]$ denotes the {\em closed interval} containing all the real numbers between $a$ and $b$.
    
    \item 
        An {\em axis-parallel box} $B$ in $\R^d$ is a set of the form $[a_1,b_1]\times[a_2,b_2]\times\dots\times[a_d,b_d]$, where $\forall i\in[d]$, $a_{i}, b_{i}\in\R$ with $a_{i} \leq b_{i}$.

    \item 
        An {\em axis-parallel $k$-flat} $K$ in $\R^d$ is a $k$-dimensional affine space such that $\exists I\subseteq [d]$ with $|I|=d-k$ and $\forall x\in K,\; \forall i\in I,\; x_i=a_i$, for some fixed number $a_i\in \RR$. Or, in other words, $K$ is the affine subspace given by $a+\text{Span}(e_{i_1},e_{i_2},\dots,e_{i_k})$, where $a\in\RR^d$ and $e_i$ is the unit vector whose $i$-th coordinate is one and the rest of the coordinates are zero.

    \item 
        A collection $\F$ of axis-parallel boxes in $\R^d$ is called a {\em $k$-collection} if $\cF$ is not pierceable by any finite family $\mathcal{K}$ of axis-parallel $k$-flats in $\R^d$.

    \item 
        A set $S \subseteq \RR^d$ is an {\em $i$-strip} if $S=S_1\times S_2\times\dots\times S_d$ and there exists $I \subseteq [d]$ with $|I| = i$ such that $\forall j\in [d]\setminus I,\; S_j=[a_j,b_j] \subseteq \RR$ with $-\infty<a_j\leq b_j<\infty$, and $\forall j\in I, S_j=\RR$.
        
        Clearly, any compact box in $\R^d$ is a {\em $0$-strip} and a {\em $d$-strip} is the whole $\R^d$. 
        We define a {\em$(-1)$-strip} to be the empty set.

    \item 
        Let $i\in[d]$. For any set $A\subseteq\R^d$, $\pi_i (A)$ denotes the projection of $A$ onto the $i$-th axis.
\end{itemize}

\section{Proof of Theorem~\ref{thm-main-theorem}}

In this section, we will prove Theorem~\ref{thm-main-theorem}. We first prove the theorem for hyperplane transversal and then prove the general case.

\begin{theorem}[$(\aleph_{0},2)$-Theorem for boxes and axis-parallel hyperplanes]
Let $\F$ be a collection of axis-parallel boxes in $\R^d$ such that $\cF$ is not pierceable by any finite number of axis-parallel hyperplanes. Then $\cF$ contains an infinite sequence $\{B_n\}_{n\in\NN}$ of boxes such that no two boxes from this sequence can be hit by an axis-parallel hyperplane.

\label{hyperplane}
\end{theorem}
\noindent
Observe that the statement of the above theorem is in contrapositive form with respect to the statement of Theorem~\ref{thm-main-theorem} (when restricted to $k = d-1$).

\begin{proof}[Proof of Theorem~\ref{hyperplane}]
Let $d-d'$ be the smallest value of $i\in\{0,1,\dots,d\}$ such that there is an $i$-strip $S$ and a $(d-1)$-collection $\F_1\subset\F$ such that $\forall B\in\cF_1,\; B\subset S$.
Without loss of generality, let a $(d-d')$-strip containing $\F_1$ be 
$$
    S_1=[a_{1,1},b_{1,1}]\times[a_{1,2},b_{1,2}]\times\dots\times[a_{1,d'},b_{1,d'}]\times\R^{d-d'}.
$$

\paragraph{Case-I:} $d'=0$. \\

\noindent In this case, $S_1=\R^d$.
Pick any $B_1\in\F$, and let 
$$
    \B_1 := \left\{ B\in\F\mid \forall i\in[d], \pi_i (B)\cap \pi_i (B_1)=\emptyset \right\}.
$$
We claim that $\B_1$ is a $(d-1)$-collection.
If not, then $\B'_1 := \F\setminus\B_1$ would be a $(d-1)$-collection.
Then, $\forall B\in\B'_1$, there is a $i \in [d]$ depending on $B$ such that $\pi_{i}(B)\cap\pi_{i}(B_1)\neq\emptyset$.
Consider the $(d-1)$-strips 
$$
    \pi_1(B_1)\times\R^{d-1},\dots,\R^{j-1}\times\pi_j(B_1)\times\R^{d-j},\dots,\R^{d-1}\times\pi_d(B_1).
$$
Since the boundaries of these $(d-1)$-strips consist of finitely many axis-parallel hyperplanes, the collection of boxes in $\B'_1$ that intersect the boundaries of these $(d-1)$-strips cannot be a $(d-1)$-collection. Let 
$\B''_1$ be the collection of all boxes in $\B'_1$ which do not intersect the boundaries of the above $(d-1)$-strips, which means,
$$
    \B''_1 \subseteq \bigcup\limits_{1\leq j\leq d} \left\{B\in\B'_1\mid B\subset\R^{j-1}\times\pi_j(B_1)\times\R^{d-j}\right\}.
$$
As $\B''_1$ is a $(d-1)$-collection, there exists a $j\in\{1,\dots,d\}$ such that
$$
    \left\{B\in\B'_1\mid B\subset\R^{j-1}\times\pi_j(B_1)\times\R^{d-j}\right\}
$$ 
is a $(d-1)$-collection. But, this would mean that there exists a $(d-1)$-strip that contains a $(d-1)$-collection which is a subset of $\F$. 
This contradicts the fact that $d' = 0$.
Therefore, $\B_1$ is a $(d-1)$-collection and we pick $B_2\in\B_1$.
Clearly, $\forall i\in[d]$, we have $\pi_i(B_2)\cap\pi_i(B_1)=\emptyset$.
Suppose we have picked $m$ boxes $B_1,B_2,\dots,B_m$ such that $B_i\in\F$ for all $i\in[m]$ and no two distinct $B_i, B_j$ can be intersected by an axis-parallel hyperplane.
Then, as argued above, we can choose a $B_{m+1}\in\F$ that does not intersect any of the finitely many $(d-1)$-strips of the form $\R^{j-1}\times\pi_j(B_i)\times\R^{d-j}$ for $j\in[d]$ and $i\in[m]$. Continuing this process we can construct an infinite sequence of boxes $\left\{ B_{n} \right\}_{n \in \mathbb{N}}$ such that for all $i \in \mathbb{N}$ we have $B_{i} \in \cF$ and
any axis-parallel hyperplane can intersect at most one box from this infinite sequence.

\paragraph{Case-II:} $d'>0$.\\

\noindent 
Let $\forall j\in[d']$, 
$$
    c_{1,j} : = \frac{b_{1,j} + a_{1,j}}{2}, \; 
    A_{1,j}^{1} : = \left[a_{1,j}, c_{1,j}\right], \; \mbox{and} \;
    A_{1,j}^{2} : = \left[ c_{1,j}, b_{1,j} \right]. 
$$
Now consider the $d'$ hyperplanes $\left\{ H_{j} \right\}_{j \in [d']}$ where $H_j = \left\{x\in\R^d\mid x_j=c_{1,j} \right\}$, $\forall j\in[d']$. Suppose $$
    \hat{\cF}=\left\{B\in\cF_1\;|\; B\cap H_j\neq\emptyset \text{ for some } j\in[d']\right\},
$$ 
that is, $\hat{\cF}$ is a collection of boxes from $\cF_1$ which are pierced by at least one hyperplane from $\left\{ H_{j} \right\}_{j \in [d']}$. As $\cF_{1}$ is a $(d-1)$-collection, $\cF_1\setminus\hat{\cF}$ must also be a $(d-1)$-collection. 
The $d'$ hyperplanes in $\left\{ H_{j}\right\}_{j \in [d']}$ split $S_1$ into at most $2^{d'}$ many $(d-d')$-strips. Therefore there exist a $(d-1)$-collection $\F_{2}$ and a $(d-d')$-strip $S_{2}$
such that 
\begin{enumerate}
    \item[(a)]
        $\F_{2} \subseteq \cF_1 \setminus \hat{\cF} \subseteq \cF_{1}$,

    \item[(b)]
        $S_2=A_{1,1}^{i_1}\times A_{1,2}^{i_2}\times\dots A_{1,d'}^{i_{d'}}\times\R^{d-d'}$ where $i_{j} \in \{1,2\}$ for all $j \in [d']$, and

    \item[(c)]
        $\forall B\in\cF_2$, we have $B\subseteq S_2$.
\end{enumerate}

Again, as $S_2$ is of the form $[a_{2,1},b_{2,1}]\times\dots\times[a_{2,d'},b_{2,d'}]\times\R^{d-d'}$ we split $S_2$, in the same way as $S_1$ was split, to obtain a $(d-d')$-strip $S_3=[a_{3,1},b_{3,1}]\times\dots\times[a_{3,d'},b_{3,d'}]\times\R^{d-d'}\subseteq S_2$, say, such that there is a $(d-1)$-collection $\F_{3} \subseteq \F_{2}$ with the property $\forall B\in\F_3$, then $B\subseteq S_{3}$.
Proceeding this way, we get that, for each $n\in\N$, there is a $(d-1)$-collection $\F_n$ and a $(d-d')$-strip $S_n=[a_{n,1},b_{n,1}]\times[a_{n,d'},b_{n,d'}]\times\R^{d-d'}$ such that the following hold:
\begin{itemize}
    \item[(i)] $\forall B\in\F_n,\; B\subseteq S_n$,
    \item[(ii)] $\F_{n+1}\subseteq \cF_{n}$,
    \item[(iii)] $S_{n+1}\subseteq S_n$ 
    \item[(iv)] $\forall i\in[d']$, $\frac{b_{n,i}-a_{n,i}}{2}=b_{n+1,i}-a_{n+1,i}$ and
    \item[(v)] $\forall i\in[d']$, $[a_{n+1,i},b_{n+1,i}]\subseteq [a_{n,i},b_{n,i}]$.
\end{itemize} 
Therefore, for each $i\in[d']$, $\exists c_i\in\RR$ such that \begin{equation}
    \lim_{n\to\infty} a_{n,i}=\lim_{n\to\infty}b_{n,i}=c_i
\label{eq: limit}\end{equation}
Now for each $i\in [d']$, consider the axis-parallel hyperplane $h_i=\{x\in\R^d\mid x_i=c_i\}$ and $\cH = \left\{ h_i\;|\; i\in [d'] \right\}$.
        
As $\F_1$ is a $(d-1)$-collection, $\cH$ is not a transversal of $\cF_{1}$. Therefore $\exists B_1\in\F_1$ such that $\forall i\in[d']$, $B_1\cap h_i=\emptyset$. Using the fact that
$\forall i\in [d']$ we have $c_i\not\in\pi_i(B_1)$ and Equation~\eqref{eq: limit}, we can show that $\exists t_1\in\N$ such that for $\forall i\in[d']$ and $\forall n\geq t_1$, 
$$
    [a_{n,i},b_{n,i}]\cap \pi_i(B_1)=\emptyset.
$$
 Again from the definition of $d'$, we know that $\forall j\in[d]\setminus[d']$, the $(d-d'-1)$-strip $S_{t_1}\cap(\R^{j-1}\times\pi_j (B_1)\times\R^{d-j})$ contains no $(d-1)$-collection that is a subset of $\F$. 
 So, we can pick a $B_2\in\F_{t_1}$ such that $\forall j\in [d]$, we have $B_2\cap (\R^{j-1}\times \pi_j(B_1)\times\R^{d-j})=\emptyset$ and $\forall h\in\cH$, we have $B_2\cap h=\emptyset$.

Continuing this process, suppose we have constructed $m$ boxes $B_{1}, \dots, B_{m}$ satisfying the following properties:
\begin{description}
    \item[{\bf Prop-1:}]
        for all $i \in [m]$ we have $B_{i} \in \cF$,

    \item[{\bf Prop-2:}] 
        no axis-parallel hyperplane pierces more than one box from the set $\left\{ B_{1}, \dots, B_{m}\right\}$, and
        
    \item[{\bf Prop-3:}] 
        $\forall s\in [m]$ and $\forall h_{i} \in \mathcal{H}$, we have
        $B_{s}\cap h_i=\emptyset$.
\end{description}
Since $\forall h_{i} \in \mathcal{H}$ and $\forall j \in [m]$, $B_{j} \cap h_{i} \neq \emptyset$, there exists a $t_{m} \in \mathbb{N}$ such that $\forall n \geq t_{m}$, $\forall r \in [m]$ and $\forall i \in [d']$, we have $[a_{n,i}\, , b_{n,i}] \cap \pi_{i}(B_{r}) = \emptyset$.
Let $\B_m$ the collection of boxes in $\F_{t_{m}}$ which do not intersect the $(d-d'-1)$-strip $S_{t_m}\cap(\R^{j-1}\times\pi_j (B_r)\times\R^{d-k})$ for any $j\in[d]\setminus [d']$ and any $r\in[m]$. Observe that $\B_m$ will be a $(d-1)$-collection. Pick $B_{m+1}\in\B_m$ so that $B_{m+1}$ does not intersect $h_i$ for any $h_{i} \in \mathcal{H}$.
Therefore we now have $m+1$ boxes $B_{1}, \dots, B_{m+1}$ that satisfy {\bf Prop-1}, {\bf Prop-2}, and {\bf Prop-3} given above with $m$ replaced by $m+1$.  Continuing this process, we obtain an infinite sequence of boxes $\left\{ B_{n} \right\}_{n \in \mathbb{N}}$ such that for all $i \in \mathbb{N}$ we have $B_{i} \in \cF$ and
any axis-parallel hyperplane can intersect at most one box from this infinite sequence. This concludes the proof of the theorem.
\end{proof}

Now we prove Theorem \ref{thm-main-theorem}.
\kflats*

\begin{proof}
        We use induction on $d$ to prove the Theorem.
        Note that we have proved Theorem \ref{thm-main-theorem} for $d=1$ and $k=0$.
        Assume that Theorem \ref{thm-main-theorem} is proven for all $d<m$ and $k<d$.
        We now prove it when $d=m$.
        The case where $k=m-1$ has already been proven, so we assume that $k<m-1$.
        \begin{case}
            \textit{There is a $k$-collection $\F'\subset\F$ which is intersected by an axis-parallel hyperplane}\\
            
            \noindent Without loss of generality, let $h= \left\{ x\in\R^m\mid x_1=a \right\}$ be an axis-parallel hyperplane that intersects all sets in a $k$-collection $\F'$ that is a subset of $\F$.
            Let $\widetilde{B}=B\cap h$ for all $B\in \F'$.
            Let $\widetilde{\F}'=\{\Tilde{B}\mid B\in\F'\}$.
            Then $\Tilde{\F'}$ is a $k$-collection in an $(m-1)$-flat, so by the induction hypothesis, there is a sequence $\{\Tilde{B_n}\}_{n=1}^{\infty}$, no two of which can be pierced by an axis-parallel $k$-flat in $h$.
            We claim that the corresponding sequence $\{B_n\}_{n=1}^{\infty}$ is the required sequence.
            If not, then let two distinct $B_s,B_j$ be pierced by an axis-parallel $k$-flat $K$ such that $x\in K$ whenever $x_{r_i}=a_i$ for all $i\in[m-k]$, where 
            $$
                1\leq r_1 < r_2 < \dots < r_{m-k}\leq m.
            $$
            If $r_1=1$, then as all sets are axis-parallel boxes, the axis-parallel $k$-flat $K'$, defined as
            $$
                K' := \left\{ x\in\R^{m} \mid x_{1}=a,\mbox{ and } x_{r_i}=a_i\mbox{ for all }i=2,\dots,m-k \right\},
            $$
            lies on $h$ and intersects $\Tilde{B_s},\Tilde{B_j}$, 
            which is a contradiction.
            If $r_1>1$, then $K\cap h$ is an axis-parallel $(k-1)$-flat on $h$ intersecting $B_s,B_j$, which is a contradiction since every axis-parallel $(k-1)$-flat is a subset of an axis-parallel $k$-flat.
            Therefore, $\{B_n\}_{n=1}^{\infty}$ is the required sequence.
        \end{case}

        \begin{case}
            \textit{There is no $k$-collection that is a subset of $\F$ which is intersected by an axis-parallel hyperplane}\\
            
            \noindent In this case, $\F$ is also a $(d-1)$-collection, because if $\F$ could be pierced by finitely many axis-parallel hyperplanes, then at least one of the hyperplanes would pierce a $k$-collection.
            Then, by Theorem~\ref{hyperplane}, we can find a sequence $\left\{ B_{n} \right\}_{n=1}^{\infty}$ in $\F$ such that no distinct pair $B_i,B_j$ with $i,j\in\N$, $i>j$ can be pierced by an axis-parallel hyperplane.
            Since an axis-parallel $k$-flat is a subset of an axis-parallel hyperplane, the sequence $\{B_n\}$ satisfies the condition of Theorem~\ref{thm-main-theorem}.
        \end{case}
\end{proof}

\section{Feasibilty of different colorful extensions}

In this section, we will investigate different colorful extensions of Theorem~\ref{thm-main-theorem}. 
First, we show the impossibility of obtaining a strictly colorful version of Theorem~\ref{thm-main-theorem}. 
\impossibility*
\begin{proof}
    Let $\QQ^d = \left\{ Q_{1}, Q_{2}, \dots \right\}$ be an ordering of the elements of $\QQ^d$ and for all $n \in \mathbb{N}$, $Q_n=\left(q_{n,1},\dots,q_{n,d}\right)$.
    For all $m$ and $n$ in $\N$, we define 
    $$
        B_{n,m} := \left[q_{n,1}-\frac{1}{2^{n+2m}},q_{n,1}-\frac{1}{2^{n+2m+1}}\right] \times \dots \times \left[q_{n,d}-\frac{1}{2^{n+2m}},q_{n,d}-\frac{1}{2^{n+2m+1}}\right],
    $$
    and 
    \begin{center}
        $\cF_{n} : = \Big\{ B_{n,m} \;\mid\; m \in \N \Big\}$ 
    \end{center}
    $$
        \cF_{n} : = \Big\{ B_{n,m} \;\mid\; m \in \N \Big\}.
    $$
    Observe that for all $i \in [d]$ and $m \neq m'$, we have 
    $$
        \pi_i(B_{n,m})\cap\pi_i(B_{n,m'})=\emptyset, \; \forall n \in \N 
    $$
    This implies that, for all $n \in \N$, $\cF_{n}$ is a $(d-1)$-collection, and therefore $\F_n$ is also a $k$-collection for every $k\in\left\{ 0, \dots, d-1 \right\}$.

    Now let, if possible, $t\in\NN$ and $\{B_n\}_{n\in\N}$ be an infinite sequence of boxes where for all $n\in\N$, $B_n\in\F_{t+n}$ and for no distinct $i,j\in\N$ there exists an axis-parallel $k$-flat that intersects both $B_i$ and $B_j$. Let $B_1 = \left[a_{1}, b_{1} \right] \times \dots \times \left[a_{d}, b_{d} \right]$, and there exists $\epsilon > 0$ such that, for all $i \in [d]$, we have $\epsilon < b_i - a_i$.
    For each $i\in[d]$, set $c_i := \frac{a_i + b_i}{2}$ and 
    $$
        B := \left[ c_1-\frac{\epsilon}{4},c_1+\frac{\epsilon}{4} \right] \times \dots \times \left[ c_d-\frac{\epsilon}{4},c_d+\frac{\epsilon}{4} \right].
    $$
    Observe that there exists $n' \in \N$ such that for all $\ell > n'$ we have 
    $$
        \frac{1}{2^{\ell}}<\frac{\epsilon}{8}.
    $$
    Since $B$ contains infinitely many elements of $\QQ^d$, there is a $Q_r = \left(q_{r,1}, \dots, q_{r,d}\right) \in B$ with $r>n'$.
    For every $m\in\N$, we have 
    $$
        B_{r,m} \subseteq \left[ q_{r,1} -\frac{1}{2^r},q_{r,1} \right] \times \dots \times \left[ q_{r,d} - \frac{1}{2^r}, q_{r,d} \right].
    $$
    Now, note that for every $i\in[d]$, we have 
    $$
        \left[ q_{r,i}-\frac{1}{2^r},q_{r,i} \right] \subseteq \left[c_i-\frac{\epsilon}{4}-\frac{\epsilon}{8},c_i+\frac{\epsilon}{4}\right] \subseteq \left[c_i-\frac{\epsilon}{2},c_i+\frac{\epsilon}{2}\right] \subseteq\left[ a_i,b_i \right].
    $$
    Therefore, for every $m\in\N$, we have $B_{r,m}\subseteq B_1$, which contradicts the claim that any axis-parallel $k$-flat can intersect at most one box from the infinite sequence $\left\{ B_{n}\right\}_{n \in \N}$.
\end{proof}
\noindent
In a later paper by Jung and P\'{a}lv\"{o}lgyi~\cite{JungP2023}, Keller offered an alternative construction that can be adapted for boxes as well which shows that the strict heterochromatic $(\aleph_0,k+2)$-property does not necessarily imply that there is a family of boxes that is finitely pierceable: let $\F_0=\{B_n\}_{n=1}^{\infty}$ be a collection of boxes that cannot be pierced by finitely many $k$-flats. 
For every $n\in\N$, let $\F_n$ be a collection of boxes that lie completely within $B_n$ and $\F_n$ is not pierceable by finitely many $k$-flats.
Then every strictly heterochromatic sequence $\{B_n\}_{n=0}^{\infty}$ chosen from $\{\F_n\}_{n=0}^{\infty}$ contains two sets that are pierced by a point even though no $\F_n$ is pierceable by finitely many points, since all sets in some $\F_n$ is contained in $B_0$.
We get Theorem~\ref{impossibility_strong_colorful} by taking infinitely many copies of $\F_0$ along with $\F_1,\F_2,\dots$ in the above construction.


However, we can prove Theorem~\ref{thm-colorful-main-theorem}.

            
        

\colorfulmaintheorem*
\begin{proof}
We shall prove the contrapositive statement, that is, if no $\cF_{\ell}$ has a finite-size axis-parallel $k$-transversal, then we get an infinite colorful sequence of boxes $\{B_n\}_{n\in\NN}$ from $\{\cF_n\}_{n\in\NN}$ such that no $3$ boxes can be pierced by an axis-parallel $k$-flat.

    We first prove Theorem~\ref{thm-colorful-main-theorem} for the case where $k=d-1$.\\
    
    \noindent \textbf{Case-I}: $k=d-1$.\\
    
    \noindent Following the proof of Theorem~\ref{hyperplane}, for every $n\in\N$ $d-d'_n$ is the smallest value of $i\in\{0,1,\dots,d\}$ such that there is an $i$-strip $S_n$ and a $(d-1)$-collection $\F'_{n} \subseteq \F_{n}$ such that $\forall B\in\F'_n$, we have $B\subseteq S_n$. We can say that there is a $r\in\{0,1,\dots,d\}$ and a set $I\subset [d]$ with $|I|=d-r$ such that for an infinite subset $J \subseteq \N$, we have:
    \begin{enumerate}
        \item 
            For every $n\in J$, we have $d'_n=r$,

        \item 
            For every $n\in J$, there is an $(d-r)$-strip $S'_n$ such that $S'_n=S'_{n,1}\times\dots\times S'_{n,d}$ with $S'_{n,i}=\R$ whenever $i\in I$, and there is a $(d-1)$-collection $\F'_{n} \subseteq \F_{n}$ with the property that $\forall B\in\F'_n$, we have $B \subseteq S'_n$.
    \end{enumerate}
    In order to reduce the burden of notations, let us assume that $J=\N$, $\F'_n=\F_n$ $\forall n\in\N$ and 
    $$
        I = \left\{ d-r+1,d-r+2,\dots,d \right\}.
    $$
    Then, following the proof of Theorem~\ref{hyperplane}, for every $n\in\N$ we get a $c_n=(c_{n,1},\dots,c_{n,r})$ such that $\forall \epsilon>0$, there is a $(d-1)$-collection $\F_{n,\epsilon}\subset\F_n$ with the property that $\forall B\in\F_{n,\epsilon}$ we have 
    $$
        B \subseteq \left[c_{n,1}-\epsilon,c_{n,1}+\epsilon \right] \times \dots \times \left[c_{n,r}-\epsilon,c_{n,r}+\epsilon \right] \times \R^{d-r}.
    $$
    Let $C = \left\{ c_{n} \right\}_{n\in\N}$.
    It is easy to see that we can obtain a subsequence $\{c_{n_m}\}_{m\in\N}$ of $C$ such that for all $i\in[r]$, either $\{c_{n_m,i}\}_{m\in\N}$ converges to some $c_i\in\R$ or diverges to either $\infty$ or $-\infty$.
    Let $I'\subseteq[r]$ be the set for which if $i\in I'$, then $c_{n_m,i}\to c_i$ as $m\to\infty$ for some $c_i\in\R$, and otherwise $c_{n_m,i}\to\infty\text{ or }c_{n_m,i}\to-\infty$.
    
    In order to avoid ugly notations, for all $m\in\N$ we write $\F_{n_m}$ as $\F_m$ and $c_{n_m}$ as $c_n$. 
    Pick $B_1\in\F_{1}$ so that $\forall i\in I'$, we have $c_i\notin\pi_i (B_1)$.
    In general, let $B_1,\dots,B_t$, $B_i\in\F_{n_i}$ be such that no two $B_i,B_j$ is intersected by an axis-parallel hyperplane, and for all $i\in I'$ and $j\in[t]$, we have $c_i\notin\pi_i (B_j)$.
    Then there is an $n_{t+1}\in\N$ and an $\epsilon>0$ for which, $\forall i\in[r] $ and $\forall j\in[t]$, we have 
    $$
        \left[ c_{n_{t+1},i}-\epsilon, c_{n_{t+1},i}+\epsilon \right]\cap\pi_i (B_j)=\emptyset.
    $$
    From the definition of $c_{n_{t+1}}$, we have that the $(d-r)$-strip 
    $$
        \left[c_{n_{t+1},1}-\epsilon,c_{n_{t+1},1}+\epsilon \right] \times \dots \times \left[c_{n_{t+1},r}-\epsilon,c_{n_{t+1},r}+\epsilon\right] \times \R^{d-r}
    $$ 
    contains a $(d-1)$-collection $\F_{n_{t+1},\epsilon}\subseteq\F_{n_{t+1}}$.
    Since $d-r$ is the smallest value of $i$ for which an $i$-strip contains a $(d-1)$-collection that is a subset of $\F_{n_{t+1}}$, we have that for every $j\in[t]$ and $i\in[d]$, there is a $(d-1)$-collection that is a subset of $\F_{n_{t+1}}$ in 
    $$
        \Big( \left[ c_{n_{t+1},1}-\epsilon,c_{n_{t+1},1}+\epsilon \right] \times\dots \times \left[c_{n_{t+1},r}-\epsilon,c_{n_{t+1},r}+\epsilon \right] \times \R^{d-r} \Big) \, \setminus \, \left( \R^{i-1}\times\pi_i (B_j)\times\R^{d-i} \right).
    $$
    Also, all sets in $\F_{n_{t+1}}$ that do not intersect the axis-parallel hyperplane $h_{i} = \left\{x\in\R^d\mid x_{i} = c_{i} \right\}$ for any $i\in I'$ form a $(d-1)$-collection, since $|I'|$ is finite.
    Therefore, we can find a $B_{t+1}\in\F_{n_{t+1}}$ such that 
    \begin{align*}
        B_{t+1} \subseteq
        \bigcap\limits_{i\in[d],j\in [t]}  \left( S \, \setminus \, \R^{i-1}\times\pi_i (B_j)\times\R^{d-i}  \right)
    \end{align*}
    where 
    $$
        S = [c_{n_{t+1},1}-\epsilon,c_{n_{t+1},1}+\epsilon]\times\dots\times[c_{n_{t+1},r}-\epsilon,c_{n_{t+1},r}+\epsilon]\times\R^{d-r},
    $$ 
    and $\forall i\in I'$ we have $c_i\notin \pi_i(B_{t+1})$.
    Proceeding this way, we obtain the required sequence $\{B_n\}_{n\in\N}$.\\

    \noindent\textbf{Case-II:} $0\leq k < d-1$.\\

    \noindent We prove this case by induction on $d$.
    We have already shown the case to hold for $d=1$ and $k=0$.
    Let the case be true for $d=1,2,\dots,D-1$ and all $k\in\{0,1,\dots,d-1\}$.
    We now prove it for $d=D$.
    If for infinitely many $n\in\N$, $\F_n$ is also a $(D-1)$-collection, then we can use the arguments of Case-I to obtain the required sequence of boxes.
    So, without loss of generality, let us assume that no $\F_n$ for any $n\in\N$ is a $(D-1)$-collection and $k<d-1$.
    Let 
    $$
        \F = \bigcup_{n\in\N} \F_{n}.
    $$
    Then either $\F$ is a $(D-1)$-collection, or it is not.
    Let us first consider the case where $\F$ is not a $(D-1)$-collection.\\
    
    \noindent\textbf{Subcase-II(a):} $\F$ is not a $(D-1)$-collection.\\
    
    \noindent Then $\F$ can be pierced by finitely many axis-parallel hyperplanes $h_1,\dots,h_t$.
    Let $\forall n\in\N$ and $i\in[t]$, $\F_{n,i}$ denote the subset of $\F_n$ such that $B\in\F_{n,i}$ if and only if $B\cap h_i\neq\emptyset$.
    Then, for each $n\in\N$, at least one $\F_{n,i}$, $i\in[t]$ is a $k$-collection, since 
    $$
        \bigcup_{i\in[t]}\F_{n,i}=\F_{n}
    $$ 
    and each $\F_n$ is a $k$-collection.
    Therefore, there is an $i\in[t]$ for which $\F_{n,i}$ is a $k$-collection for infinitely many $n\in\N$. Without loss of generality, let $i=1$ and let $\{n_m\}_{m\in\N}$ be a strictly increasing sequence in $\N$ such that $\F_{n_m,1}$ is a $k$-collection. 
    In order to avoid ugly notations, we write $\F_{n_m,1}$ as $\F_m$ for each $m\in\N$.
    Let 
    $$
        \widetilde{\F}_{n} = \left\{ B\cap h_{1} \mid B \in \F_{n} \right\}
    $$ 
    and for every $B\in\F_n$, let 
    $$
        \widetilde{B} = B \cap h_{1}.
    $$
    Then $\widetilde{\F}_{n}$ is a $k$-collection in $\R^{D-1}$, and therefore by the induction hypothesis, there is a sequence $\{n_m\}_{m\in\N}$ and $\{B_m\}_{m\in\N}$ such that $B_m\in\F_{n_m}$ and there is no axis-parallel $k$-flat on $h_1$ that pierces $\widetilde{B}_m$ and $\widetilde{B}_{m'}$ when $m\neq m'$.
    As argued in the proof of Theorem~\ref{thm-main-theorem}, we can show that no two $B_m,B_{m'}$, $m\neq m'$ can be hit by an axis-parallel $k$-flat in $\R^D$.\\

    \noindent\textbf{Subcase-II(b):} $\F$ is a $(D-1)$-collection.\\
    
    \noindent In this case, by Theorem~\ref{thm-main-theorem} we get a sequence $\{B_m\}_{m\in\N}$ such that no two $B_m,B_{m'}$ can be hit by an axis-parallel hyperplane, which also means that $B_m,B_{m'}$ cannot be hit by any axis-parallel $k$-flat.
    We show that the sequence $\{B_m\}_{m\in\N}$ comes from infinitely many $\F_n$'s.
    Define, for all $n\in\N$, 
    $$
        \mathcal{G}_{n} : = \bigcup_{i\in[n]}\F_{i}.
    $$
    Then, as $\mathcal{G}_{n}$ is the union of finitely many $k$-collections, none of which is a $(d-1)$-collection, we have that $\mathcal{G}_n$ is not a $(d-1)$-collection for any $n\in\N$.
    But as $\{B_m\}_{m\in\N}$ is a $(d-1)$-collection, $\{B_m\}_{m\in\N}$ cannot be a subset of $\mathcal{G}_n$ for any $n\in\N$.
    This shows that $\{B_m\}_{m\in\N}$ comes from infinitely many $\F_n$'s.\\ \\
    \noindent This concludes the proof.
\end{proof}

\appendix

\section{Arguments from Chakraborty et al.~\cite{chakraborty2024heterochromatic}}
\label{sec:appendix1}

Here we provide the arguments from Chakraborty, Ghosh, and Nandi~\cite{chakraborty2024heterochromatic} that shows that the existence of a heterochromatic $(\aleph_0,q)$-theorem implies the existence of a finite hitting set which works for all but finitely many families.
Let $\cB$ and $\cS$ be collections of nonempty sets in $\RR^d$.
For any $\cF\subseteq\cB$, we say $\cF$ is {\em finitely pierceable by $\cS$} if there exists a finite set $S\subseteq \cS$ such that for all $B\in\cF$, we have $\bigcup_{A\in S} A\cap B \neq \emptyset$.
A sequence $\{\cF_n\}_{n\in\NN}$ with $\cF_n\subseteq\cB$ $\forall n\in\NN$ is said to have the {\em heterochromatic $(\aleph_0,q)$-property with respect to $\mathcal{S}$} if for every heterochromatic sequence $\{B_n\}_{n\in\NN}$ chosen from $\{\cF_n\}_{n\in\NN}$, we have a set $A\in\cS$ such that there are $q$ distinct $B_n$'s for which $B_n\cap A\neq \emptyset$.
We say that the {\em heterochromatic $(\aleph_0,q)$-theorem for $\cB$ with respect to $\cS$ holds} if for any $\{\cF_n\}_{n\in\NN}$ with $\cF_n\subseteq\cB$ satisfying the heterochromatic $(\aleph_0,q)$-property with respect to $\mathcal{S}$, there exists an $n\in\NN$ for which $\cF_n$ is finitely pierceable by $\cS$.
We say that the {\em infinite heterochromatic $(\aleph_0,q)$-theorem for $\cB$ with respect to $\cS$ holds} if for any $\{\cF_n\}_{n\in\NN}$ with $\cF_n\subseteq\cB$ satisfying the heterochromatic $(\aleph_0,q)$-property with respect to $\mathcal{S}$, all but finitely many $\cF_n$'s are finitely pierceable by $\cS$.

\begin{lemma}[Heterochromatic Property Lemma]
    The heterochromatic $(\aleph_0,q)$-theorem for $\cB$ with respect to $\cS$ implies the infinite heterochromatic $(\aleph_0,q)$-theorem for $\cB$ with respect to $\cS$.
    \label{lem: heterochromatic}
\end{lemma}

\begin{proof}
Let the heterochromatic $(\aleph_0,q)$-theorem for $\cB$ with respect to $\cS$ hold.
Let $\{\cF_n\}_{n\in\NN}$, $\cF_n\subseteq\cB$, be a sequence that satisfies the heterochromatic $(\aleph_0,q)$-property with respect to $\mathcal{S}$.
If possible, let $\{\cF_{n_i}\}_{i\in\NN}$ where $\{n_i\}_{i\in\NN}$ is a strictly increasing sequence in $\NN$ be such that no $\cF_{n_i}$ is finitely pierceable by $\cS$.
As $\{\cF_{n_i}\}_{i\in\NN}$ also satisfies the heterochromatic $(\aleph_0,q)$-property with respect to $\mathcal{S}$, there is an $i\in\NN$ for which $\cF_{n_i}$ is finitely pierceable by $\cS$.
But this is a contradiction, and therefore, the result follows.
\end{proof}
So, we shall be using the phrases {\em heterochromatic $(\aleph_0,q)$-theorem} and {\em infinite heterochromatic $(\aleph_0,q)$-theorem} interchangeably.
It is trivially seen that for some $\cB$ and $\cS$ satisfying the heterochromatic $(\aleph_0,q)$-theorem and any $N\in\NN$, we can find a sequence $\{\cF_n\}_{n\in\NN}$ with $\cF_n\subseteq\cB$ that satisfy the heterochromatic $(\aleph_0,q)$ property and for $N$ distinct values of $n\in\NN$, $\cF_{n}$ is not finitely pierceable by $\cS$.
A simple example is found by setting $\cB$ to be the set of all unit balls in $\RR^d$, $\cS$ to be the set of all points in $\RR^d$ and $\cF_1=\dots=\cF_N=\cB$ and $\cF_n=\{B(\mathcal{O},1)\}$ for all $n\in\NN\setminus\{1,\dots,N\}$.

Now we show that the heterochromatic $(\aleph_0,k+2)$-theorem implies the existence of a finite hitting set that works for all but finitely many families.
\begin{lemma}[Heterochromatic Theorem and a Finite-size Piercing Set]\label{lem: hetero_means_finite}
Let $\cB$ and $\cS$ be collections of sets in $\R^d$ and $q\in\NN$, where $\cB$ satisfies the heterochromatic $(\aleph_0,q)$-theorem with respect to $\cS$.
If $\{\F_n\}_{n\in\NN}$ be a sequence of subsets of $\cB$ that satisfies the heterochromatic $(\aleph_0,q)$-property with respect to $\cS$, then there exists an $N\in\NN$ and a finite collection $\cK\subset\cS$ such that $\forall n>N$ where $n\in\NN$, $\forall C\in\F_n$, we have $\bigcup_{K\in\cK}K\cap C\neq\emptyset$.
\end{lemma}

\begin{proof}
By Lemma~\ref{lem: heterochromatic} we have that there is an $N\in\NN$ such that every $\F_n$ with $n>N$ is finitely pierceable by $\cS$. 
Let $\F=\bigcup_{n=N}^{\infty}\F_n$.
Then for any infinite set $F\subseteq\F$, either $F$ contains infinitely many elements of some $\F_n$ with $n\geq N$, or $F$ contains a heterochromatic sequence with respect to the sequence $\left\{ \F_n \right\}_{n=N}^{\infty}$.
In either case, this means there is an $S\in\cS$ that pierces $q$ elements of $F$.
Since the heterochromatic $(\aleph_0,q)$-theorem implies the monochromatic $(\aleph_0,q)$-theorem, $\F$ must be finitely pierceable by $\cS$.
\end{proof}

\bibliographystyle{alpha}

\bibliography{references}

\newcommand{\etalchar}[1]{$^{#1}$}
\begin{thebibliography}{DLGMM19}

\bibitem[ADLS17]{AmentaDeLS17}
N.~Amenta, J.~A De~Loera, and P.~Sober{\'o}n.
\newblock {Helly’s Theorem: New Variations and Applications}.
\newblock {\em {Algebraic and Geometric Methods in Discrete Mathematics}},
  685:55--95, 2017.

\bibitem[AGP02]{AGP}
B.~Aronov, J.~E. Goodman, and R.~Pollack.
\newblock {A Helly-Type Theorem for Higher-Dimensional Transversals}.
\newblock {\em Computational Geometry}, 21(3):177--183, 2002.

\bibitem[AGPW00]{ABJP}
B.~Aronov, J.~E. Goodman, R.~Pollack, and R.~Wenger.
\newblock {A Helly-Type Theorem for Hyperplane Transversals to Well-Separated
  Convex Sets}.
\newblock In {\em Proceedings of the 16th Annual Symposium on Computational
  Geometry, SoCG}, page 57–63, 2000.

\bibitem[AK92]{AlonK92}
N.~Alon and D.~J. Kleitman.
\newblock {Piercing Convex Sets and the Hadwiger-Debrunner $(p,q)$-Problem}.
\newblock {\em Advances in Mathematics}, 96(1):103--112, 1992.

\bibitem[AK95]{AlonK95}
N.~Alon and G.~Kalai.
\newblock {Bounding the Piercing Number}.
\newblock {\em Discrete \& Computational Geometry}, 13:245--256, 1995.

\bibitem[AKMM02]{AKMM}
N.~Alon, G.~Kalai, J.~Matoušek, and R.~Meshulam.
\newblock Transversal numbers for hypergraphs arising in geometry.
\newblock {\em Advances in Applied Mathematics}, 29(1):79--101, 2002.

\bibitem[Alo98]{Alon98}
N.~Alon.
\newblock {Piercing $d$-Intervals}.
\newblock {\em Discrete \& Computational Geometry}, 19(3):333 -- 334, 1998.

\bibitem[Alo12]{Alon12}
N.~Alon.
\newblock {A Non-linear Lower Bound for Planar Epsilon-nets}.
\newblock {\em Discrete \& Computational Geometry}, 47(2):235 -- 244, 2012.

\bibitem[B{\'{a}}r82]{Barany82}
I.~B{\'{a}}r{\'{a}}ny.
\newblock {A Generalization of Carath{\'{e}}odory's Theorem}.
\newblock {\em Discrete Mathematics}, 40(2-3):141 -- 152, 1982.

\bibitem[BK22]{BaranyK2022helly}
I.~B{\'a}r{\'a}ny and G.~Kalai.
\newblock {Helly-type Problems}.
\newblock {\em Bulletin of the American Mathematical Society}, 59(4):471--502,
  2022.

\bibitem[CD20]{ChenD20}
Ke~Chen and Adrian Dumitrescu.
\newblock {On Wegner's inequality for axis-parallel rectangles}.
\newblock {\em Discrete Mathematics}, 343(12):112091, 2020.

\bibitem[CFPS15]{CorreaFPS15}
J.~R. Correa, L.~Feuilloley, P.~P{\'{e}}rez{-}Lantero, and J.~A. Soto.
\newblock {Independent and Hitting Sets of Rectangles Intersecting a Diagonal
  Line: Algorithms and Complexity}.
\newblock {\em Discrete \& Computational Geometry}, 53(2):344--365, 2015.

\bibitem[CGHP08]{CheongGHP08}
O.~Cheong, X.~Goaoc, A.~F. Holmsen, and S.~Petitjean.
\newblock {Helly-Type Theorems for Line Transversals to Disjoint Unit Balls}.
\newblock {\em Discrete \& Computational Geometry}, 39(1-3):194 -- 212, 2008.

\bibitem[CGN23]{chakraborty2023stabbing}
S.~Chakraborty, A.~Ghosh, and S.~Nandi.
\newblock {Stabbing boxes with finitely many axis-parallel lines and flats}.
\newblock {\em CoRR}, abs/2308.10479, 2023.

\bibitem[CGN24a]{chakraborty2024countably}
S.~Chakraborty, A.~Ghosh, and S.~Nandi.
\newblock {Countably Colorful Hyperplane Transversal}.
\newblock {\em CoRR}, abs/2402.10012, 2024.

\bibitem[CGN24b]{chakraborty2024heterochromatic}
S.~Chakraborty, A.~Ghosh, and S.~Nandi.
\newblock {Heterochromatic Higher Order Transversals for Convex Sets}.
\newblock {\em CoRR}, abs/2212.14091, 2024.

\bibitem[CH21]{ChanH21}
T.~M. Chan and S.~Har{-}Peled.
\newblock {Smallest $k$-Enclosing Rectangle Revisited}.
\newblock {\em Discrete \& Computational Geometry}, 66(2):769--791, 2021.

\bibitem[CL01]{ChazelleL01}
B.~Chazelle and A.~Lvov.
\newblock {The Discrepancy of Boxes in Higher Dimension}.
\newblock {\em Discrete \& Computational Geometry}, 25(4):519 -- 524, 2001.

\bibitem[CPST09]{ChenPST09}
X.~Chen, J.~Pach, M.~Szegedy, and G.~Tardos.
\newblock {Delaunay graphs of point sets in the plane with respect to
  axis-parallel rectangles}.
\newblock {\em Random Structures \& Algorithms}, 34(1):11 -- 23, 2009.

\bibitem[CSZ18]{ChudnovskySZ18}
M.~Chudnovsky, S.~Spirkl, and S.~Zerbib.
\newblock {Piercing Axis-Parallel Boxes}.
\newblock {\em The Electronic Journal of Combinatorics}, 25(1):1, 2018.

\bibitem[CW21]{ChalermsookW21}
P.~Chalermsook and B.~Walczak.
\newblock {Coloring and Maximum Weight Independent Set of Rectangles}.
\newblock In {\em Proceedings of the 2021 {ACM-SIAM} Symposium on Discrete
  Algorithms, {SODA}}, pages 860 -- 868, 2021.

\bibitem[DG82]{DanzerG82}
L.~Danzer and B.~Gr{\"{u}}nbaum.
\newblock {Intersection Properties of Boxes in $\mathbb{R}^{d}$}.
\newblock {\em Combinatorica}, 2(3):237--246, 1982.

\bibitem[DGK63]{DGK}
L.~Danzer, B.~Gr{\"u}nbaum, and V.~Klee.
\newblock {\em Helly's Theorem and Its Relatives}.
\newblock Proceedings of Symposia in Pure Mathematics: Convexity. American
  Mathematical Society, 1963.

\bibitem[DLGMM19]{de2019discrete}
J.~De~Loera, X.~Goaoc, F.~Meunier, and N.~Mustafa.
\newblock {The discrete yet ubiquitous theorems of Carath{\'e}odory, Helly,
  Sperner, Tucker, and Tverberg}.
\newblock {\em Bulletin of the American Mathematical Society}, 56(3):415--511,
  2019.

\bibitem[DT15]{DumitrescuT15}
A.~Dumitrescu and C.~D. T{\'{o}}th.
\newblock Packing anchored rectangles.
\newblock {\em Combinatorica}, 35(1):39--61, 2015.

\bibitem[FK93]{Fon-Der-FlaassK93}
D.~Fon{-}Der{-}Flaass and A.~V. Kostochka.
\newblock {Covering boxes by points}.
\newblock {\em Discrete Mathematics}, 120(1-3):269--275, 1993.

\bibitem[GKM{\etalchar{+}}22]{GalvezKMMPW22}
W.~G{\'{a}}lvez, A.~Khan, M.~Mari, T.~M{\"{o}}mke, M.~R. Pittu, and A.~Wiese.
\newblock {A $3$-Approximation Algorithm for Maximum Independent Set of
  Rectangles}.
\newblock In {\em Proceedings of the 2022 {ACM-SIAM} Symposium on Discrete
  Algorithms, {SODA}}, pages 894 -- 905, 2022.

\bibitem[GL85]{GyarfasL85}
A.~Gy{\'{a}}rf{\'{a}}s and J.~Lehel.
\newblock {Covering and coloring problems for relatives of intervals}.
\newblock {\em Discrete Mathematics}, 55(2):167--180, 1985.

\bibitem[GP88]{GoodmanPJAMS88}
J.~Goodman and R.~Pollack.
\newblock {Hadwiger's Transversal Theorem In Higher Dimensions}.
\newblock {\em Journal of The American Mathematical Society}, 1:301 -- 301,
  1988.

\bibitem[Had56]{Hadwiger56}
H.~Hadwiger.
\newblock {\"{U}ber einen Satz Hellyscher Art}.
\newblock {\em Arch. Math.}, 7:377--379, 1956.

\bibitem[Had57]{Hadwiger57}
H.~Hadwiger.
\newblock \"{U}ber eibereiche mit gemeinsamer treffgeraden.
\newblock {\em Portugal Math.}, 6:23 -- 29, 1957.

\bibitem[HD57]{HadwigerD57}
H.~Hadwiger and H.~Debrunner.
\newblock {\"{U}ber eine Variante zum Hellyschen Satz}.
\newblock {\em Archiv der Mathematik}, 8(4):309--313, 1957.

\bibitem[Hel23]{Helly23}
E.~Helly.
\newblock {\"{U}ber Mengen konvexer K\"{o}rper mit gemeinschaftlichen Punkten}.
\newblock {\em Jahresbericht der Deutschen Mathematiker-Vereinigung}, 32:175 --
  176, 1923.

\bibitem[HKL03]{HolmsenKL03}
A.~F. Holmsen, M.~Katchalski, and T.~Lewis.
\newblock {A Helly-Type Theorem for Line Transversals to Disjoint Unit Balls}.
\newblock {\em Discrete \& Computational Geometry}, 29(4):595 -- 602, 2003.

\bibitem[HW17]{HolmsenW17}
A.~Holmsen and R.~Wenger.
\newblock {Helly-type Theorems and Geometric Transversals}.
\newblock In J.~O'Rourke, J.~E. Goodman, and C.~D. T\'{o}th, editors, {\em
  {Handbook of Discrete and Computational Geometry, 3rd Edition}}, chapter~4,
  pages 91--123. CRC Press LLC, Boca Raton, FL, 2017.

\bibitem[JP23]{JungP2023}
A.~Jung and D.~P{\'{a}}lv{\"{o}}lgyi.
\newblock $k$-dimensional transversals for balls.
\newblock {\em CoRR}, abs/2311.15646, 2023.

\bibitem[KMP16]{KupavskiiMP16}
A.~Kupavskii, N.~H. Mustafa, and J.~Pach.
\newblock {New Lower Bounds for epsilon-Nets}.
\newblock In {\em Proceedings of the 32nd International Symposium on
  Computational Geometry, SoCG}, volume~51, pages 54:1 -- 54:16, 2016.

\bibitem[KNS03]{KatzNS03}
M.~J. Katz, F.~Nielsen, and M.~Segal.
\newblock {Maintenance of a Piercing Set for Intervals with Applications}.
\newblock {\em Algorithmica}, 36(1):59--73, 2003.

\bibitem[KP22]{KellerP22}
C.~Keller and M.~A. Perles.
\newblock {An $(\aleph_{0}, k+2)$-Theorem for $k$-Transversals}.
\newblock In {\em Proceedings of the 38th International Symposium on
  Computational Geometry, {SoCG}}, volume 224, pages 50:1--50:14, 2022.

\bibitem[KP23]{KellerP23}
C.~Keller and M.~A. Perles.
\newblock {An $(\aleph_0,k+2)$-Theorem for $k$-Transversals}.
\newblock {\em CoRR}, abs/2306.02181, 2023.

\bibitem[KS20]{KellerS20}
C.~Keller and S.~Smorodinsky.
\newblock {From a $(p, 2)$-Theorem to a Tight $(p, q)$-Theorem}.
\newblock {\em Discrete \& Computational Geometry}, 63(4):821 -- 847, 2020.

\bibitem[MDG18]{MustafaDG18}
N.~H. Mustafa, K.~Dutta, and A.~Ghosh.
\newblock {A Simple Proof of Optimal Epsilon Nets}.
\newblock {\em Combinatorica}, 38(5):1269 -- 1277, 2018.

\bibitem[Mit21]{Mitchell21}
J.~S.~B. Mitchell.
\newblock {Approximating Maximum Independent Set for Rectangles in the Plane}.
\newblock In {\em Proceedings of the 62nd {IEEE} Annual Symposium on
  Foundations of Computer Science, {FOCS}}, pages 339 -- 350, 2021.

\bibitem[MR17]{MustafaR17}
N.~H. Mustafa and S.~Ray.
\newblock { $\varepsilon$-Mnets: Hitting Geometric Set Systems with Subsets}.
\newblock {\em Discrete \& Computational Geometry}, 57(3):625 -- 640, 2017.

\bibitem[Nie00]{Nielsen00}
F.~Nielsen.
\newblock {Fast stabbing of boxes in high dimensions}.
\newblock {\em Theoretical Computer Science}, 246(1-2):53--72, 2000.

\bibitem[PT10]{PachT10}
J.~Pach and G.~Tardos.
\newblock Coloring axis-parallel rectangles.
\newblock {\em Journal of Combinatorial Theory, Series A}, 117(6):776 -- 782,
  2010.

\bibitem[PT12]{PachT12}
J.~Pach and G.~Tardos.
\newblock {Piercing quasi-rectangles - On a problem of Danzer and Rogers}.
\newblock {\em Journal of Combinatorial Theory, Series A}, 119(7):1391 -- 1397,
  2012.

\bibitem[PT13]{PachT13}
J.~Pach and G.~Tardos.
\newblock {Tight lower bounds for the size of epsilon-nets}.
\newblock {\em Journal of the American Mathematical Society}, 26:645 -- 658,
  2013.

\bibitem[RST10]{RafalinST10}
E.~Rafalin, D.~L. Souvaine, and C.~D. T{\'{o}}th.
\newblock {Cuttings for Disks and Axis-Aligned Rectangles in Three-Space}.
\newblock {\em Discrete \& Computational Geometry}, 43(2):221 -- 241, 2010.

\bibitem[San40]{Santal2009}
L~Santal{\'o}.
\newblock Un teorema sobre conjuntos de paralelepipedos de aristas paralelas.
\newblock {\em Publ. Inst. Mat. Univ. Nac. Litoral}, 2:49--60, 1940.

\bibitem[Tom23]{tomon2022lower}
I.~Tomon.
\newblock Lower bounds for piercing and coloring boxes.
\newblock {\em Advances in Mathematics}, 435:109360, 2023.

\bibitem[T{\'{o}}t08]{Toth08}
C.~D. T{\'{o}}th.
\newblock {Binary Space Partitions for Axis-Aligned Fat Rectangles}.
\newblock {\em SIAM Journal on Computing}, 38(1):429 -- 447, 2008.

\end{thebibliography}

\end{document}